# MODELING THE NON-LINEAR VISCOELASTIC RESPONSE OF HIGH TEMPERATURE POLYIMIDES

SATISH KARRA AND K. R. RAJAGOPAL

Abstract. A constitutive model is developed to predict the viscoelastic response of polyimide resins that are used in high temperature applications. This model is based on a thermodynamic framework that uses the notion that the 'natural configuration' of a body evolves as the body undergoes a process and the evolution is determined by maximizing the rate of entropy production in general and the rate of dissipation within purely mechanical considerations. We constitutively prescribe forms for the specific Helmholtz potential and the rate of dissipation (which is the product of density, temperature and the rate of entropy production), and the model is derived by maximizing the rate of dissipation with the constraint of incompressibility, and the reduced energy dissipation equation is also regarded as a constraint in that it is required to be met in every process that the body undergoes. The efficacy of the model is ascertained by comparing the predictions of the model with the experimental data for PMR-15 and HFPE-II-52 polyimide resins.

## 1. Introduction

Polyimides are well known to be extremely stable at high temperatures and also have a glass transition temperature that is greater than 300°C. Due to their good performance at high temperature ranges they are used by aircraft and automobile industries to fashion their products. They are also used in wafer fabrication due to their excellent high temperature resistance and adhesive properties (see Ghosh and Mittal (1996)). The mechanical properties of polyimides and polyimide composites used in several applications especially in the aerospace industry are affected by high temperature, diffusion of moisture and subsequent oxidation. Hence, there is need for a good understanding of the various degradation mechanisms that are operational when such materials are subject to hostile environment. Recent experimental evidence shows that the response of polyimide resins is solid-like viscoelastic response (see Bhargava (2007), Falcone and Ruggles-Wrenn (2009)). Moreover, the response of such bodies is non-linear (For a detailed description concerning how to differentiate between viscoelastic solid-like and fluid-like response we refer the reader to Wineman and Rajagopal (2000), Málek and Rajagopal (2006)). A thermodynamic framework which takes into account the viscoelastic solid-like response of the polyimides along with the various degradation processes needs to be developed. As a first step, our aim in this paper is to develop a model based on a thermodynamic framework that can predict the non-linear viscoelastic solid-like response of polyimides at various temperatures. This first step is non-trivial and presents interesting challenges. Subsequently, we shall extend this model to include degradation due to moisture diffusion, and chemical reactions, specifically oxidation.

While the thermal response of linear viscoelastic solids have been studied in great detail, there has been no systematic study of non-linear viscoelastic solids. Standard techniques like superposition that are valid in linear response are no longer valid, thus making the study much more complicated. A single integral model has been proposed by Pipkin and Rogers (1968) who have assumed that a linear combination of responses to single step strain histories can be used as an approximation to the response to an arbitrary strain history. Unfortunately, such a model does not have a sound thermodynamic basis, and moreover the model is too general to be of use. Later on, Fung (1993) developed a quasi-linear viscoelastic model that has been shown to predict the behavior of several biological materials, though not adequately when the strains are large. This model by Fung can be shown to be a special case of the model by Pipkin and Rodgers. For further details on the various viscoelastic models for solids that have been reported in the literature, see the review articles by Drapaca et al. (2007) and Wineman (2009).

Most of the literature concerning the modeling of the response of polyimides use a viscoelastic model proposed by Schapery (1969). Muliana and Sawant (2009) used Schapery's model and obtained material parameters for PMR-15 [1] using the experiments carried out by Marais and Villoutreix (1998). They used these material parameters in their 'micromechanical' model to predict the behavior of Kevlar/PMR-15 composites. Ahci and Talreja (2006) performed experiments on a composite made of graphite fiber in a HFPE-II[2]. They have extended the framework developed by Schapery to include a 'damage tensor' as an internal variable, and have also included anisotropy to model the composite behavior. Recently, Falcone and Ruggles-Wrenn performed experiments on PMR-15 at the service temperature of an aircraft, 288°C, and compared the predictions of Schapery's model for the problem of creep with experimental data. Bhargava (2007) has also used Schapery's model to predict the behavior of HFPE-II-52. More recently, Hall (2008) developed a thermodynamic framework for finite anisotropic viscoplastic models to study the response of polymers subject to extreme thermal environment.

Given its extensive use some comments on Schapery's model are warranted: Schapery developed a viscoelastic model using linear phenomenological relations based on Onsager's reciprocity theorem (see Onsager (1931)) which states that the forces are linearly related to the fluxes near equilibrium (see for instance, equation 11 in Schapery (1964)). Next, he introduced nonlinearity by assuming that the coefficient matrix relating the forces and fluxes depends on generalized coordinates and temperature. Furthermore, the free energy expression was obtained using a Taylor series expansion and by neglecting higher order terms. Thus, while the model might be able to describe slight deviation from linear response, one cannot expect it to be capable of describing truly non-linear response that is thermodynamically compatible. Thus, if one is interested in describing the non-linear response of viscoelastic solids that takes into account its thermodynamic effects, a different model is necessary. In this paper we develop a viscoelastic solid model based on a thermodynamic framework that can be used to describe the non-linear response exhibited by a class of polymers. The framework has been recently developed and is used to describe the response of bodies that produce entropy in a variety of ways. In order to derive meaningful physical models they require that amongst the class of processes that are possible the process which is actually taken by the body is one that maximizes the rate of entropy production. One can find details concerning this approach in the review article by

[1] PMR polymerization of monomer reactant.

[2] HFPE stands for hydrofluoropolyether.

Rajagopal and Srinivasa (2004a). In this approach, one need not assume near-equilibrium behavior and linear phenomenological relations between forces and fluxes, the approach is much more general. Also, one need not use a Taylor series expansion of the free energy, and neglect higher order terms. Recently, Rajagopal and Srinivasa (2004b) have shown that if one uses an expression for entropy production which is quadratic in the fluxes, one can arrive at Onsager's relations upon maximizing the rate of entropy production along with appropriate constraints. As mentioned earlier such a thermodynamic framework has also been used to model various material responses such as viscoelastic solid-like and fluid-like behavior, traditional plasticity, twinning, crystallization and so on (see the review article by Rajagopal and Srinivasa (2004a) for the references and for the details of the framework).

In this paper, a viscoelastic solid model is derived by assuming forms for the Helmholtz potential and the rate of dissipation, and maximizing the rate of dissipation [3] with incompressibility and the reduced energy dissipation equation as constraints. This model is shown to predict the viscoelastic response of polyimide resin. Experimental data for PMR-15 polyimide resin from Falcone and Ruggles-Wrenn (2009), and for HFPE-II-52 from Bhargava (2007) are used to evaluate the efficacy of the model.

The current paper is organized as follows. In section (2), the kinematics that is required in this paper are documented. In sections (3.1), (3.2), a viscoelastic solid model is developed using a thermodynamic framework. We show that the viscoelastic solid model that is developed is a generalization of the one-dimensional standard linear solid model in (3.3). In section (3.4) the problem of uniaxial extension is set up using our model, and the creep solution obtained by using the model is compared with experimental data for PMR-15 and HFPE-II-52 polyimide resins in section (3.5). We find that the theoretical predictions agree quite well with the experimental results.

## 2. Preliminaries

Let $\kappa_R(\mathcal{B})$ and $\kappa_t(\mathcal{B})$ denote the reference configuration and the current configuration, respectively. The motion $\chi_{\kappa_R}$ is defined as the one-one mapping that assigns to each point $\boldsymbol{X} \in \kappa_R$, a point $\boldsymbol{x} \in \kappa_t$, at a time $t$, i.e.,

$$\boldsymbol{x} = \chi_{\kappa_R}(\boldsymbol{X}, t). \tag{2.1}$$

The mapping $\chi_{\kappa_R}(\boldsymbol{X}, t)$ is assumed to be sufficiently smooth and invertible. Let $\kappa_{p(t)}$ be the stress-free configuration instantaneously reached by the body upon removal of the external stimuli. We assume that the body can be instantaneously unloaded. We shall call this configuration as the natural configuration corresponding to $\kappa_t$. The natural configuration that underlies the current configuration depends on the process class that is admissible. Thus underlying natural configuration corresponding to isothermal and adiabatic processes could be different. Let $\boldsymbol{F}$ be gradient of motion $\chi_{\kappa_R}(\boldsymbol{X}, t)$ (usually known as the deformation gradient), defined by

$$\boldsymbol{F} := \frac{\partial \chi_{\kappa_R}}{\partial \boldsymbol{X}}, \tag{2.2}$$

and let the left and right Cauchy-Green tensors be defined through

$$\boldsymbol{B} = \boldsymbol{F}\boldsymbol{F}^T, \quad \boldsymbol{C} = \boldsymbol{F}^T\boldsymbol{F}. \tag{2.3}$$

[3]In case of isothermal processes, the rate of dissipation is the rate of conversion of mechanical working into heat (energy in thermal form), and in general it is the product of density, temperature and the rate of entropy production.

Let $\boldsymbol{F}_{\kappa_{p(t)}}$ be the gradient of the mapping from $\kappa_{p(t)}$ to $\kappa_t$, and let $\boldsymbol{G}$ be defined by

$$\boldsymbol{G} := \boldsymbol{F}_{\kappa_R \to \kappa_{p(t)}} = \boldsymbol{F}^{-1}_{\kappa_{p(t)}} \boldsymbol{F}. \tag{2.4}$$

Similar to (2.3), we shall denote the left Cauchy-Green stretch tensors $\boldsymbol{B}_G$ and $\boldsymbol{B}_{p(t)}$ as

$$\boldsymbol{B}_G := \boldsymbol{G}\boldsymbol{G}^T, \quad \boldsymbol{B}_{p(t)} := \boldsymbol{F}_{\kappa_{p(t)}} \boldsymbol{F}^T_{\kappa_{p(t)}}. \tag{2.5}$$

We shall also define the velocity gradients

$$\boldsymbol{L}_G := \dot{\boldsymbol{G}}\boldsymbol{G}^{-1}, \quad \boldsymbol{L} = \dot{\boldsymbol{F}}\boldsymbol{F}^{-1}, \quad \boldsymbol{L}_p = \dot{\boldsymbol{F}}_{\kappa_{p(t)}} \boldsymbol{F}^{-1}_{\kappa_{p(t)}}, \tag{2.6}$$

and their symmetric parts by

$$\boldsymbol{D}_i = \frac{1}{2}\left(\boldsymbol{L}_i + \boldsymbol{L}_i^T\right), \quad i = p(t), G \quad \text{or no subscript}. \tag{2.7}$$

Also, we define the principal invariants through

$$\mathrm{I}_{\boldsymbol{B}_l} = \mathrm{tr}(\boldsymbol{B}_l), \quad \mathrm{II}_{\boldsymbol{B}_l} = \frac{1}{2}\left\{[\mathrm{tr}(\boldsymbol{B}_l)]^2 - \mathrm{tr}(\boldsymbol{B}_l^2)\right\}, \quad \mathrm{III}_{\boldsymbol{B}_l} = \det(\boldsymbol{B}_l) \quad l = G, p(t), \tag{2.8}$$

where tr(.) is the trace operator for a second order tensor and det(.) is the determinant. Now, from (2.4)

$$\begin{aligned} \dot{\boldsymbol{F}} &= \dot{\boldsymbol{F}}_{\kappa_{p(t)}} \boldsymbol{G} + \boldsymbol{F}_{\kappa_{p(t)}} \dot{\boldsymbol{G}} \\ \Rightarrow \dot{\boldsymbol{F}}\boldsymbol{F}^{-1} &= \dot{\boldsymbol{F}}_{\kappa_{p(t)}} \boldsymbol{G}\boldsymbol{G}^{-1}\boldsymbol{F}^{-1}_{\kappa_{p(t)}} + \boldsymbol{F}_{\kappa_{p(t)}} \dot{\boldsymbol{G}} \\ \Rightarrow \boldsymbol{L} &= \boldsymbol{L}_{p(t)} + \boldsymbol{F}_{\kappa_{p(t)}} \boldsymbol{L}_G \boldsymbol{F}^{-1}_{\kappa_{p(t)}}, \end{aligned} \tag{2.9}$$

where $\dot{(.)}$ is the material time derivative of the second order tensor. In addition,

$$\begin{aligned} \dot{\boldsymbol{B}}_{p(t)} &= \dot{\boldsymbol{F}}_{\kappa_{p(t)}} \boldsymbol{F}^T + \boldsymbol{F}\dot{\boldsymbol{F}}^T_{\kappa_{p(t)}} \\ &= \boldsymbol{L}_{p(t)} \boldsymbol{B}_{p(t)} + \boldsymbol{B}_{p(t)} \boldsymbol{L}^T_{p(t)}, \end{aligned} \tag{2.10}$$

and similarly

$$\dot{\boldsymbol{B}}_G = \boldsymbol{L}_G \boldsymbol{B}_G + \boldsymbol{B}_G \boldsymbol{L}_G^T. \tag{2.11}$$

Hence, from (2.9) and (2.10), we have

$$\dot{\boldsymbol{B}}_{p(t)} = \boldsymbol{L}\boldsymbol{B}_{p(t)} + \boldsymbol{B}_{p(t)} \boldsymbol{L}^T_{p(t)} - \boldsymbol{F}_{\kappa_{p(t)}} \left(\boldsymbol{L}_G + \boldsymbol{L}_G^T\right) \boldsymbol{F}^T_{\kappa_{p(t)}}, \tag{2.12}$$

and so

$$\overset{\triangledown}{\boldsymbol{B}}_{p(t)} = -2\boldsymbol{F}_{\kappa_{p(t)}} \boldsymbol{D}_G \boldsymbol{F}^T_{\kappa_{p(t)}}, \tag{2.13}$$

where $\overset{\triangledown}{(.)}$ is the usual Oldroyd derivative defined through $\overset{\triangledown}{\boldsymbol{A}} := \dot{\boldsymbol{A}} - \boldsymbol{L}\boldsymbol{A} - \boldsymbol{A}\boldsymbol{L}^T$. When one considers non-isothermal processes the local form of the second law of thermodynamics takes the following form:

$$\boldsymbol{T} \cdot \boldsymbol{D} - \varrho\dot{\psi} - \varrho s\dot{\theta} - \frac{\boldsymbol{q}_h \cdot \mathrm{grad}(\theta)}{\theta} = \varrho\theta\zeta := \xi \geq 0, \tag{2.14}$$

where $\boldsymbol{T}$ is the Cauchy stress, $\psi$ is the specific Helmholtz potential, $\varrho$ is the density, $\theta$ is the temperature, $s$ is the specific entropy, $\boldsymbol{q}_h$ is the heat flux, $\zeta$ is the rate of entropy production and $\xi$ is the rate of dissipation.

## 3. Constitutive assumptions and maximization of the rate of dissipation

3.1. **General results.** We shall assume that the viscoelastic solid is isotropic and incompressible with the specific Helmholtz potential of the form

$$\psi = \psi(\boldsymbol{B}_{p(t)}, \boldsymbol{B}_G, \theta) = \hat{\psi}(\mathrm{I}_{\boldsymbol{B}_{p(t)}}, \mathrm{II}_{\boldsymbol{B}_{p(t)}}, \mathrm{I}_{\boldsymbol{B}_G}, \mathrm{II}_{\boldsymbol{B}_G}, \theta). \tag{3.1}$$

Since the elastic response is isotropic, without loss of generality, we choose $\kappa_{p(t)}$ such that

$$\boldsymbol{F}_{\kappa_{p(t)}} = \boldsymbol{V}_{\kappa_{p(t)}}, \tag{3.2}$$

where $\boldsymbol{V}_{\kappa_{p(t)}}$ is the right stretch tensor in the polar decomposition of $\boldsymbol{F}_{\kappa_{p(t)}}$. We shall also assume that the total rate of dissipation can be split additively as follows

$$\boldsymbol{T} \cdot \boldsymbol{D} - \varrho\dot{\psi} - \varrho s\dot{\theta} = \xi_m \geq 0, \quad -\frac{\boldsymbol{q}_h \cdot \mathrm{grad}(\theta)}{\theta} = \xi_c \geq 0, \tag{3.3}$$

where $\xi_m$, $\xi_c$ are the rates of mechanical dissipation (conversion of working into thermal energy) and dissipation due to heat conduction, respectively. Now, we constitutively choose

$$\boldsymbol{q}_h = -k(\theta)\,\mathrm{grad}(\theta), \quad k(\theta) \geq 0, \tag{3.4}$$

where $k$ is the thermal conductivity, so that $3.3_{(b)}$ is automatically satisfied.

Next,

$$\begin{aligned}\dot{\psi} = &\left[\left(\frac{\partial\hat{\psi}}{\partial \mathrm{I}_{B_{p(t)}}} + \mathrm{I}_{B_{p(t)}}\frac{\partial\hat{\psi}}{\partial \mathrm{II}_{B_{p(t)}}}\right)\boldsymbol{I} - \frac{\partial\hat{\psi}}{\partial \mathrm{II}_{B_{p(t)}}}\boldsymbol{B}_{p(t)}\right] \cdot \dot{\boldsymbol{B}}_{p(t)} \\ &+ \left[\left(\frac{\partial\hat{\psi}}{\partial \mathrm{I}_{B_G}} + \mathrm{I}_{B_G}\frac{\partial\hat{\psi}}{\partial \mathrm{II}_{B_G}}\right)\boldsymbol{I} - \frac{\partial\hat{\psi}}{\partial \mathrm{II}_{B_G}}\boldsymbol{B}_G\right] \cdot \dot{\boldsymbol{B}}_G + \frac{\partial\hat{\psi}}{\partial\theta}\dot{\theta},\end{aligned} \tag{3.5}$$

and using (2.10), (2.11) along with (3.2) in (3.5), we obtain

$$\begin{aligned}\dot{\psi} = 2&\left[\left(\frac{\partial\hat{\psi}}{\partial \mathrm{I}_{B_{p(t)}}} + \mathrm{I}_{B_{p(t)}}\frac{\partial\hat{\psi}}{\partial \mathrm{II}_{B_{p(t)}}}\right)\boldsymbol{B}_{p(t)} - \frac{\partial\hat{\psi}}{\partial \mathrm{II}_{B_{p(t)}}}\boldsymbol{B}^2_{p(t)}\right] \cdot (\boldsymbol{D} - \boldsymbol{D}_G) \\ + 2&\left[\left(\frac{\partial\hat{\psi}}{\partial \mathrm{I}_{B_G}} + \mathrm{I}_{B_G}\frac{\partial\hat{\psi}}{\partial \mathrm{II}_{B_G}}\right)\boldsymbol{B}_G - \frac{\partial\hat{\psi}}{\partial \mathrm{II}_{B_G}}\boldsymbol{B}^2_G\right] \cdot \boldsymbol{D}_G + \frac{\partial\hat{\psi}}{\partial\theta}\dot{\theta}.\end{aligned} \tag{3.6}$$

Next, we shall assume the rate of mechanical dissipation to be of the form

$$\xi_m = \xi_m(\theta, \boldsymbol{B}_{p(t)}, \boldsymbol{D}_G). \tag{3.7}$$

On substituting (3.6) into ($3.3_{(a)}$), we arrive at

$$\begin{aligned}
&\left[\boldsymbol{T} - 2\varrho\left(\frac{\partial\hat{\psi}}{\partial \mathrm{I}_{B_{p(t)}}} + \mathrm{I}_{B_{p(t)}}\frac{\partial\hat{\psi}}{\partial \mathrm{II}_{B_{p(t)}}}\right)\boldsymbol{B}_{p(t)} + 2\varrho\frac{\partial\hat{\psi}}{\partial \mathrm{II}_{B_{p(t)}}}\boldsymbol{B}^2_{p(t)}\right]\cdot\boldsymbol{D}\\
&+ 2\varrho\left[\left(\frac{\partial\hat{\psi}}{\partial \mathrm{I}_{B_{p(t)}}} + \mathrm{I}_{B_{p(t)}}\frac{\partial\hat{\psi}}{\partial \mathrm{II}_{B_{p(t)}}}\right)\boldsymbol{B}_{p(t)} - \frac{\partial\hat{\psi}}{\partial \mathrm{II}_{B_{p(t)}}}\boldsymbol{B}^2_{p(t)}\right]\cdot\boldsymbol{D}_G\\
&- 2\varrho\left[\left(\frac{\partial\hat{\psi}}{\partial \mathrm{I}_{B_G}} + \mathrm{I}_{B_G}\frac{\partial\hat{\psi}}{\partial \mathrm{II}_{B_G}}\right)\boldsymbol{B}_G - \frac{\partial\hat{\psi}}{\partial \mathrm{II}_{B_G}}\boldsymbol{B}^2_G\right]\cdot\boldsymbol{D}_G\\
&- \varrho\left[\frac{\partial\hat{\psi}}{\partial\theta} + s\right]\dot{\theta}\\
&= \xi_m(\theta, \boldsymbol{B}_{p(t)}, \boldsymbol{D}_G).
\end{aligned} \tag{3.8}$$

We shall set

$$s = -\frac{\partial\hat{\psi}}{\partial\theta}, \tag{3.9}$$

and define

$$\boldsymbol{T}_{p(t)} := 2\varrho\left[\left(\frac{\partial\hat{\psi}}{\partial \mathrm{I}_{B_{p(t)}}} + \mathrm{I}_{B_{p(t)}}\frac{\partial\hat{\psi}}{\partial \mathrm{II}_{B_{p(t)}}}\right)\boldsymbol{B}_{p(t)} - \frac{\partial\hat{\psi}}{\partial \mathrm{II}_{B_{p(t)}}}\boldsymbol{B}^2_{p(t)}\right], \tag{3.10}$$

$$\boldsymbol{T}_G := 2\varrho\left[\left(\frac{\partial\hat{\psi}}{\partial \mathrm{I}_{B_G}} + \mathrm{I}_{B_G}\frac{\partial\hat{\psi}}{\partial \mathrm{II}_{B_G}}\right)\boldsymbol{B}_G - \frac{\partial\hat{\psi}}{\partial \mathrm{II}_{B_G}}\boldsymbol{B}^2_G\right]. \tag{3.11}$$

Using (3.9)–(3.11) in (3.8), we obtain

$$\begin{aligned}
&\big(\boldsymbol{T} - \boldsymbol{T}_{p(t)}\big)\cdot\boldsymbol{D} + \big(\boldsymbol{T}_{p(t)} - \boldsymbol{T}_G\big)\cdot\boldsymbol{D}_G\\
&\quad = \xi_m(\theta, \boldsymbol{B}_{p(t)}, \boldsymbol{D}_G).
\end{aligned} \tag{3.12}$$

From constraint of incompressibility, we have

$$\mathrm{tr}(\boldsymbol{D}) = \mathrm{tr}(\boldsymbol{D}_{p(t)}) = \mathrm{tr}(\boldsymbol{D}_G) = 0. \tag{3.13}$$

Since, RHS of (3.12) does not depend on $\boldsymbol{D}$, using (3.13),

$$\boldsymbol{T} = p\boldsymbol{I} + \boldsymbol{T}_{p(t)}, \tag{3.14}$$

where $p$ is the Lagrange multiplier due to the constraint of incompressibility, with

$$\big(\boldsymbol{T}_{p(t)} - \boldsymbol{T}_G\big)\cdot\boldsymbol{D}_G = \xi_m(\theta, \boldsymbol{B}_{p(t)}, \boldsymbol{D}_G), \tag{3.15}$$

which can be re-written as

$$(\boldsymbol{T} - \boldsymbol{T}_G)\cdot\boldsymbol{D}_G = \xi_m(\theta, \boldsymbol{B}_{p(t)}, \boldsymbol{D}_G), \tag{3.16}$$

using (3.13) and (3.14) .

Now, we shall maximize the rate of dissipation $\xi_m$ by varying $\boldsymbol{D}_G$ for fixed $\boldsymbol{B}_{p(t)}$. That is, we maximize the function[4]

$$\Phi := \xi_m + \lambda_1\left[\xi_m - (\boldsymbol{T} - \boldsymbol{T}_G)\cdot\boldsymbol{D}_G\right] + \lambda_2(\boldsymbol{I}\cdot\boldsymbol{D}_G), \tag{3.17}$$

[4]Though we only document that the first derivative is zero here, it can be shown that the extremum is a maximum.

where $\lambda_1, \lambda_2$ are the Lagrange multipliers. By setting, $\partial\Phi/\partial \boldsymbol{D}_G = 0$, we get

$$\boldsymbol{T} = \boldsymbol{T}_G + \frac{\lambda_2}{\lambda_1}\boldsymbol{I} + \left(\frac{\lambda_1 + 1}{\lambda_1}\right)\frac{\partial \xi_m}{\partial \boldsymbol{D}_G}. \tag{3.18}$$

We need to determine the Lagrange multipliers. On substituting (3.18) into (3.16), we get

$$\left(\frac{\lambda_1 + 1}{\lambda_1}\right) = \frac{\xi_m}{\frac{\partial \xi_m}{\partial \boldsymbol{D}_G}\cdot \boldsymbol{D}_G}, \tag{3.19}$$

and so (3.18) with (3.11) becomes

$$\boldsymbol{T} = 2\varrho\left[\left(\frac{\partial\hat{\psi}}{\partial \mathrm{I}_{B_G}} + \mathrm{I}_{B_G}\frac{\partial\hat{\psi}}{\partial \mathrm{II}_{B_G}}\right)\boldsymbol{B}_G - \frac{\partial\hat{\psi}}{\partial \mathrm{II}_{B_G}}\boldsymbol{B}_G^2\right] + \left(\frac{\xi_m}{\frac{\partial \xi_m}{\partial \boldsymbol{D}_G}\cdot \boldsymbol{D}_G}\right)\frac{\partial \xi_m}{\partial \boldsymbol{D}_G} + \hat{\lambda}\boldsymbol{I}, \tag{3.20}$$

where $\hat{\lambda} := \frac{\lambda_2}{\lambda_1}$ is the Lagrange multiplier due to the constraint of incompressibility.

Finally, the constitutive relations for the viscoelastic solid are given by

$$\boldsymbol{T} = p\boldsymbol{I} + 2\varrho\left[\left(\frac{\partial\hat{\psi}}{\partial \mathrm{I}_{B_{p(t)}}} + \mathrm{I}_{B_{p(t)}}\frac{\partial\hat{\psi}}{\partial \mathrm{II}_{B_{p(t)}}}\right)\boldsymbol{B}_{p(t)} - \frac{\partial\hat{\psi}}{\partial \mathrm{II}_{B_{p(t)}}}\boldsymbol{B}_{p(t)}^2\right], \tag{3.21a}$$

$$\boldsymbol{T} = \hat{\lambda}\boldsymbol{I} + 2\varrho\left[\left(\frac{\partial\hat{\psi}}{\partial \mathrm{I}_{B_G}} + \mathrm{I}_{B_G}\frac{\partial\hat{\psi}}{\partial \mathrm{II}_{B_G}}\right)\boldsymbol{B}_G - \frac{\partial\hat{\psi}}{\partial \mathrm{II}_{B_G}}\boldsymbol{B}_G^2\right] + \left(\frac{\xi_m}{\frac{\partial \xi_m}{\partial \boldsymbol{D}_G}\cdot \boldsymbol{D}_G}\right)\frac{\partial \xi_m}{\partial \boldsymbol{D}_G}, \tag{3.21b}$$

$$\boldsymbol{q}_h = -k(\theta)\mathrm{grad}(\theta), \quad s = -\frac{\partial\hat{\psi}}{\partial\theta}. \tag{3.21c}$$

3.2. **Specific case.** Specifically, we choose the specific Helmholtz potential as

$$\hat{\psi} = A^s + (B^s + c_2^s)(\theta - \theta_s) - \frac{c_1^s}{2}(\theta - \theta_s)^2 - c_2^s\theta \ln\left(\frac{\theta}{\theta_s}\right) + \frac{\mu_{G0} - \mu_{G1}\theta}{2\varrho\theta_s}(\mathrm{I}_{B_G} - 3) + \frac{\mu_{p0} - \mu_{p1}\theta}{2\varrho\theta_s}(\mathrm{I}_{B_{p(t)}} - 3), \tag{3.22}$$

where $\mu_G, \mu_p$ are elastic constants, $\theta_s$ is a reference temperature for the viscoelastic solid, and the rate of dissipation as

$$\xi_m = \eta(\theta)\left(\boldsymbol{D}_G \cdot \boldsymbol{B}_{p(t)}\boldsymbol{D}_G\right), \tag{3.23}$$

where $\eta$ is the viscosity.

Now,

$$\begin{aligned} s &= -\frac{\partial\hat{\psi}}{\partial\theta} \\ &= -(B^s + c_2^s) + c_1^s(\theta - \theta_s) + c_2^s \ln\left(\frac{\theta}{\theta_s}\right) + c_2^s + \frac{\mu_{G1}}{2\varrho\theta_s}(\mathrm{I}_{B_G} - 3) + \frac{\mu_{p1}}{2\varrho\theta_s}(\mathrm{I}_{B_{p(t)}} - 3). \end{aligned} \tag{3.24}$$

The internal energy $\epsilon$ is given by

$$\begin{aligned} \epsilon &= \hat{\psi} + \theta s \\ &= A^s - B^s\theta_s + c_2^s(\theta - \theta_s) + \frac{c_1^s}{2}\left(\theta^2 - \theta_s^2\right) + \frac{\mu_{G0}}{2\varrho\theta_s}(\mathrm{I}_{B_G} - 3) + \frac{\mu_{p0}}{2\varrho\theta_s}(\mathrm{I}_{B_{p(t)}} - 3). \end{aligned} \tag{3.25}$$

and the specific heat capacity $C_v$ is

$$C_v = \frac{\partial\epsilon}{\partial\theta} = c_1^s\theta + c_2^s. \tag{3.26}$$

Also, (3.21a), (3.21b) reduce to

$$\boldsymbol{T} = p\boldsymbol{I} + \bar{\mu}_p \boldsymbol{B}_{p(t)}, \tag{3.27a}$$

$$\boldsymbol{T} = \lambda\boldsymbol{I} + \bar{\mu}_G \boldsymbol{B}_G + \frac{\eta}{2}\left(\boldsymbol{B}_{p(t)}\boldsymbol{D}_G + \boldsymbol{D}_G\boldsymbol{B}_{p(t)}\right), \tag{3.27b}$$

where $\bar{\mu}_p = \frac{\mu_{p0}-\mu_{p1}\theta}{\theta_s}$, $\bar{\mu}_G = \frac{\mu_{G0}-\mu_{G1}\theta}{\theta_s}$. From (3.27)

$$(p-\lambda)\boldsymbol{I} + \bar{\mu}_p \boldsymbol{B}_{p(t)} = \bar{\mu}_G \boldsymbol{B}_G + \frac{\eta}{2}\left(\boldsymbol{B}_{p(t)}\boldsymbol{D}_G + \boldsymbol{D}_G\boldsymbol{B}_{p(t)}\right), \tag{3.28}$$

and so by pre-multiplying the above equation by $\boldsymbol{B}_{p(t)}^{-1}$ and taking the trace, we get

$$(p-\lambda) = \frac{\bar{\mu}_G \mathrm{tr}(\boldsymbol{B}_{p(t)}^{-1}\boldsymbol{B}_G) - 3\bar{\mu}_p}{\mathrm{tr}(\boldsymbol{B}_{p(t)}^{-1})}. \tag{3.29}$$

Using (3.29) in (3.28), we arrive at the following equation that holds:

$$\left[\frac{\bar{\mu}_G \mathrm{tr}(\boldsymbol{B}_{p(t)}^{-1}\boldsymbol{B}_G) - 3\bar{\mu}_p}{\mathrm{tr}(\boldsymbol{B}_{p(t)}^{-1})}\right]\boldsymbol{I} + \bar{\mu}_p \boldsymbol{B}_{p(t)} = \bar{\mu}_G \boldsymbol{B}_G + \frac{\eta}{2}\left(\boldsymbol{B}_{p(t)}\boldsymbol{D}_G + \boldsymbol{D}_G\boldsymbol{B}_{p(t)}\right), \tag{3.30}$$

which can be re-written as

$$\left[\frac{\bar{\mu}_G \mathrm{tr}(\boldsymbol{B}_{p(t)}^{-1}\boldsymbol{B}_G) - 3\bar{\mu}_p}{\mathrm{tr}(\boldsymbol{B}_{p(t)}^{-1})}\right]\boldsymbol{I} + \bar{\mu}_p \boldsymbol{B}_{p(t)} = \bar{\mu}_G \boldsymbol{B}_G - \frac{\eta}{4}\left(\boldsymbol{V}_{p(t)} \overset{\triangledown}{\boldsymbol{B}}_{p(t)} \boldsymbol{V}_{\kappa_{p(t)}}^{-1} + \boldsymbol{V}_{\kappa_{p(t)}}^{-1} \overset{\triangledown}{\boldsymbol{B}}_{p(t)} \boldsymbol{V}_{p(t)}\right), \tag{3.31}$$

where we have used (2.13) and (3.30). Thus, with the current choice of the specific Helmholtz potential and the rate of dissipation, we arrive at the following constitutive equations:

$$\boldsymbol{T} = p\boldsymbol{I} + \bar{\mu}_p \boldsymbol{B}_{p(t)}, \tag{3.32}$$

where the evolution of the natural configuration is given by (3.31). Also, note that the above model reduces to the generalized Maxwell fluid model derived by Rajagopal and Srinivasa (2000) when $\bar{\mu}_G = 0$. This is interesting, but not totally surprising, that we obtain a fluid model by eliminating a energy storage mechanism. In the corresponding one dimensional model this is tantamount to a spring being removed.

3.3. **Relationship to the Standard Linear Solid.** Now, (3.32), (3.31) can be re-written as

$$\boldsymbol{T} = (p + \bar{\mu}_p)\boldsymbol{I} + \bar{\mu}_p(\boldsymbol{B}_{p(t)} - \boldsymbol{I}), \tag{3.33a}$$

$$\left[\frac{\bar{\mu}_G \mathrm{tr}(\boldsymbol{B}_{p(t)}^{-1}\boldsymbol{B}_G) - 3\bar{\mu}_p}{\mathrm{tr}(\boldsymbol{B}_{p(t)}^{-1})} + \bar{\mu}_p - \bar{\mu}_G\right]\boldsymbol{I} + \bar{\mu}_p(\boldsymbol{B}_{p(t)} - \boldsymbol{I}) = \bar{\mu}_G(\boldsymbol{B}_G - \boldsymbol{I}) - \frac{\eta}{4}\left(\boldsymbol{V}_{p(t)} \overset{\triangledown}{\boldsymbol{B}}_{p(t)} \boldsymbol{V}_{\kappa_{p(t)}}^{-1} + \boldsymbol{V}_{\kappa_{p(t)}}^{-1} \overset{\triangledown}{\boldsymbol{B}}_{p(t)} \boldsymbol{V}_{p(t)}\right). \tag{3.33b}$$

If $\lambda_i$ $(i = G, p)$ is the one-dimensional stretch and $\varepsilon_i = \ln\lambda_i$ $(i = G, p)$ is the logarithmic strain, when one is restricted to one-dimension, (3.33) reduces to

$$\sigma = \bar{\mu}_p(\lambda_p^2 - 1), \tag{3.34a}$$

$$\bar{\mu}_p(\lambda_p^2 - 1) = \bar{\mu}_G(\lambda_G^2 - 1) + \eta\lambda_p^2\frac{\dot{\lambda}_G}{\lambda_G}, \tag{3.34b}$$

where $\sigma$ is the one dimensional stress. Equation (3.34) under the assumption that $\varepsilon_i \ll 1$ $(i = G, p)$ reduces to

$$\sigma = 2\bar{\mu}_p\varepsilon_p, \tag{3.35a}$$

$$2\bar{\mu}_p\varepsilon_p = 2\bar{\mu}_G\varepsilon_G + \hat{\eta}\dot{\varepsilon}_G, \tag{3.35b}$$

where $\hat{\eta} = \eta\lambda_p^2$ is the stretch dependent viscosity. Equation (3.35) can also be obtained by using a Kelvin-Voigt element (with spring constant $2\bar{\mu}_G$, viscosity of $\hat{\eta}$) and a spring (of spring constant $2\bar{\mu}_p$) in series, which is the spring-dashpot analogy for the standard linear solid. However, in this model the viscosity is stretch dependent and hence the model is a generalization of the classical standard linear solid as the viscosity in the standard linear solid model is assumed to be a constant. Hence, the viscoelastic solid model given by (3.32), (3.31) is a three-dimensional generalization of the standard linear solid. Of course, there can be infinity of three dimensional generalizations of a one dimensional model (see Karra and Rajagopal (2009)). Recently, Kannan and Rajagopal (2004) have also derived a three-dimensional viscoelastic solid model, that is different from the model developed in this paper, that also reduces to the standard linear solid. That more than one, in fact, infinity of generalizations are possible is akin to the situation in elementary mathematics and stems from the fact that infinity of three dimensional functions can have the same one dimensional projection. In fact, even when one considers the thermodynamical formulation that is used in the paper, using different forms for the specific Helmholtz potential and the rate of dissipation, and by maximizing the rate of dissipation with the necessary constraints more than one three-dimensional model reduces to the same one-dimensional model (see Karra and Rajagopal (2009) for details of an example).

3.4. **Application of the model.** Let us study the uniaxial extension, given by

$$x = \lambda(t)X, \quad y = \frac{1}{\sqrt{\lambda(t)}}Y, \quad z = \frac{1}{\sqrt{\lambda(t)}}Z, \tag{3.36}$$

within the context of this model. The velocity gradient is given by

$$\boldsymbol{L} = \text{diag}\left\{\frac{\dot{\lambda}}{\lambda}, -\frac{\dot{\lambda}}{2\lambda}, -\frac{\dot{\lambda}}{2\lambda}\right\}. \tag{3.37}$$

We shall assume that the stretch $\boldsymbol{B}_{p(t)}$ is given by

$$\boldsymbol{B}_{p(t)} = \text{diag}\left\{B, \frac{1}{\sqrt{B}}, \frac{1}{\sqrt{B}}\right\}. \tag{3.38}$$

So,

$$\dot{\boldsymbol{B}}_{p(t)} = \text{diag}\left\{\dot{B}, -\frac{\dot{B}}{2B^{3/2}}, -\frac{\dot{B}}{2B^{3/2}}\right\}, \tag{3.39}$$

$$\overset{\triangledown}{\boldsymbol{B}}_{p(t)} = \text{diag}\left\{\dot{B} - \frac{2B\dot{\lambda}}{\lambda}, -\frac{\dot{B}}{2B^{3/2}} + \frac{\dot{\lambda}}{\lambda\sqrt{B}}, -\frac{\dot{B}}{2B^{3/2}} + \frac{\dot{\lambda}}{\lambda\sqrt{B}}\right\}, \tag{3.40}$$

$$\boldsymbol{V}_{\kappa_{p(t)}} = \text{diag}\left\{\sqrt{B}, \frac{1}{B^{1/4}}, \frac{1}{B^{1/4}}\right\}, \tag{3.41}$$

and

$$\boldsymbol{D}_G = -\frac{1}{2}\text{diag}\left\{\frac{\dot{B}}{B} - \frac{2\dot{\lambda}}{\lambda}, \frac{\dot{\lambda}}{\lambda} - \frac{\dot{B}}{2B}, \frac{\dot{\lambda}}{\lambda} - \frac{\dot{B}}{2B}\right\}. \tag{3.42}$$

Also,

$$\begin{aligned} \boldsymbol{G} &= \boldsymbol{V}^{-1}_{\kappa_{p(t)}}\boldsymbol{F} \\ &= \text{diag}\left\{\frac{\lambda}{\sqrt{B}}, \frac{B^{1/4}}{\sqrt{\lambda}}, \frac{B^{1/4}}{\sqrt{\lambda}}\right\}, \end{aligned} \tag{3.43}$$

which yields

$$\boldsymbol{B}_G = \text{diag}\left\{\frac{\lambda^2}{B}, \frac{\sqrt{B}}{\lambda}, \frac{\sqrt{B}}{\lambda}\right\}. \tag{3.44}$$

and

$$\boldsymbol{B}^{-1}_{p(t)}\boldsymbol{B}_G = \text{diag}\left\{\frac{\lambda^2}{B^2}, \frac{B}{\lambda}, \frac{B}{\lambda}\right\}. \tag{3.45}$$

Substituting (3.38), (3.40), (3.45) into (3.30)

$$\frac{\dot{B}}{2} = \frac{B\dot{\lambda}}{\lambda} + \frac{\bar{\mu}_G}{\eta}\frac{\lambda^2}{B} - \frac{\bar{\mu}_p}{\eta}B - \left\{\frac{\frac{\bar{\mu}_G}{\eta}(\lambda^3 + 2B^3) - 3\frac{\bar{\mu}_p}{\eta}\lambda B^2}{\lambda B\,(1 + 2B^{3/2})}\right\}, \tag{3.46}$$

which can be re-written in the following form:

$$\dot{\lambda} = \lambda\left\{\frac{\dot{B}}{2B} - \left[\frac{1}{\eta B}\left(\bar{\mu}_G\frac{\lambda^2}{B} - \bar{\mu}_p B - \left(\frac{\bar{\mu}_G(\lambda^3 + 2B^3) - 3\bar{\mu}_p B^2\lambda}{B\lambda(1 + 2B^{3/2})}\right)\right)\right]\right\} \tag{3.47}$$

Now, from (3.38), (3.32), and using the fact that lateral surfaces are traction free, we conclude that

$$T_{11} = \bar{\mu}_p\left(B - \frac{1}{\sqrt{B}}\right). \tag{3.48}$$

We shall also use logarithmic strain (or true strain) $\varepsilon = \ln\lambda$ as our strain measure in what follows.

3.5. **Comparison with experimental creep data.** For the loading process, with known constant applied stress $T_{11}$ and material properties, (3.48) was first solved for $B(t)$. Then, (3.46) was solved with the initial condition $\lambda(0) = \sqrt{B(0)}$. For the unloading process, $T_{11}$ was set to zero and $B(t)$ was evaluated using (3.48). Then, using $\lambda(t_u^+) = \frac{\lambda(t_u^-)}{\lambda(0)}$ as the initial condition (where $t_u$ is the time when unloading starts), $\lambda(t)$ during the unloading process is evaluated using (3.46). All the ODEs were solved in MATLAB using the `ode45` solver.

In order to obtain the material parameters for a given set of experimental creep data, `fminsearch` function in MATLAB (which uses Nedler-Mead simplex method) was used to minimize the error defined by

$$error = w \times \sqrt{\frac{\sum \left(\varepsilon_{theo,load} - \varepsilon_{exp,load}\right)^2}{\sum \left(\varepsilon_{exp,load}\right)^2}} + (1-w) \times \sqrt{\frac{\sum \left(\varepsilon_{theo,unload} - \varepsilon_{exp,unload}\right)^2}{\sum \left(\varepsilon_{exp,unload}\right)^2}}, \qquad (3.49)$$

where $\varepsilon_{theo}$ denotes the theoretical strain values, $\varepsilon_{exp}$ denotes the experimental strain values, the suffixes *load*, *unload* denote the values during loading and unloading processes respectively, $w$ is a weight. The material parameters for the model was obtained for HFPE-II-52 polyimide resin using the experimental creep data from Bhargava (2007) at different temperatures (285°C, 300°C, 315°C and 330°C). To determine the efficacy of the model the following process was followed. At 285°C, the experimental data values for the loading of 0.45 UTS were used to obtain the material parameters by minimizing the error (3.49). Then, these material parameters were used for the model prediction at the other loadings of 0.30 UTS and 0.15 UTS. The loading values corresponding to the sets of experimental data which were used to obtain the material parameters at the other temperatures are shown in table (1). Similar to the process described above for 285°C, the material parameters shown in table (1) were used to predict the creep at other temperatures. The model predictions compare well with the experimental data as shown in figures (2), (3).

Next, the creep solution that stems from our model is compared to the experimental creep data of Falcone and Ruggles-Wrenn (2009) for PMR-15 resin at 288°C in figure (4). The best-fit values of the parameters for were found to be $\bar{\mu}_G = 4.42 \times 10^8$ Pa, $\bar{\mu}_p = 3.76 \times 10^8$ Pa, $\eta = 6.22 \times 10^{12}$ Pa.s. A weight of $w = 0.75$ was used since there are fewer data points for the unloading process. As it can be seen from figure (4), our model shows a good fit with the experiment. However, there is no additional experiment with which the predictive capability of the model can be tested.

It is seen in the experiments that at a loading close to the failure values, the experimental data shows permanent set in the body, and there seems to be 'yielding'. Our model being a viscoelastic solid model it cannot predict such a permanent set. Thus, the model should be generalized to take into account the inelastic response of the polymer, but this is a daunting problem that requires a careful and separate study.

In conclusion, a viscoelastic model has been developed which predicts the behavior of polyimide resins, that takes into account the thermal response, quite well. Our work can be extended to include anisotropy, and inelasticity to predict the response of polyimide composites and one can include various degradation mechanisms as well.

## Acknowledgements

The authors thank AFOSR for supporting this work through contract no: FA 9550-04-1-0137. Part of this work was done when Satish Karra was appointed as a lecturer during his Phd by the Department of Mechanical Engineering at Texas A&M University. He appreciates this support by the department.

## References

Ahci, E., Talreja, R., 2006. Characterization of viscoelasticity and damage in high temperature polymer matrix composites. Compos. Sci. Technol. 66 (14), 2506–2519.

Bhargava, P., 2007. High temperature properties of HFPE-II-52 polyimide resin and composites. Ph.D. thesis, Cornell University.

Drapaca, C. S., Sivaloganathan, S., Tenti, G., 2007. Nonlinear constitutive laws in viscoelasticity. Math. Mech. Solids 12 (5), 475–501.

Falcone, C. M., Ruggles-Wrenn, M. B., 2009. Rate dependence and short-term creep behavior of a thermoset polymer at elevated temperature. J. Pressure Vessel Technol. 131, 011403(1–8).

Fung, Y. C., 1993. Biomechanics: Mechanical Properties of Living Tissues. Springer-Verlag, New York.

Ghosh, M. K., Mittal, K. L., 1996. Polyimides: Fundamentals and Applications. Marcel Deckker, Inc., New York.

Hall, R. B., 2008. Combined thermodynamics approach for anisotropic, finite deformation overstress models of viscoplasticity. Int. J. Eng. Sci. 46 (2), 119–130.

Kannan, K., Rajagopal, K. R., 2004. A thermomechanical framework for the transition of a viscoelastic liquid to a viscoelastic solid. Math. Mech. Solids 9 (1), 37.

Karra, S., Rajagopal, K. R., 2009. Development of three dimensional constitutive theories based on lower dimensional experimental data. Appl. Math. 54 (2), 147–176.

Málek, J., Rajagopal, K. R., 2006. Mathematical issues concerning the Navier-Stokes equations and some of its generalizations. Handbook of Differential Equations: Evolutionary Equations 2, 371–459.

Marais, C., Villoutreix, G., 1998. Analysis and modeling of the creep behavior of the thermostable PMR-15 polyimide. J. Appl. Polym. Sci. 69 (10), 1983–1991.

Muliana, A. H., Sawant, S., 2009. Responses of viscoelastic polymer composites with temperature and time dependent constituents. Acta Mech. 204 (3), 155–173.

Onsager, L., 1931. Reciprocal Relations in Irreversible Processes. I. Phys. Rev. 37 (4), 405–426.

Pipkin, A. C., Rogers, T. G., 1968. A non-linear integral representation for viscoelastic behaviour. J. Mech. Phys. Solids 16 (1), 59–72.

Rajagopal, K. R., Srinivasa, A. R., 2000. A thermodynamic frame work for rate type fluid models. J. Non-Newtonian Fluid Mech. 88 (3), 207–227.

Rajagopal, K. R., Srinivasa, A. R., 2004a. On the thermomechanics of materials that have multiple natural configurations part I: viscoelasticity and classical plasticity. Z. Angew. Math. Phys. 55 (5), 861–893.

Rajagopal, K. R., Srinivasa, A. R., 2004b. On thermomechanical restrictions of continua. Proc. R. Soc. A 460 (2042), 631–651.

Schapery, R. A., 1964. Application of thermodynamics to thermomechanical, fracture, and birefringent phenomena in viscoelastic media. J. Appl. Phys. 35, 1451–1465.

Schapery, R. A., 1969. On the characterization of nonlinear viscoelastic materials. Polym. Eng. Sci. 9 (4), 295–310.

Wineman, A. S., 2009. Nonlinear viscoelastic solids–a review. Math. Mech. Solids 14 (3), 300–366.

Wineman, A. S., Rajagopal, K. R., 2000. Mechanical Response of Polymers: an Introduction. Cambridge University Press, New York.

| Temperature | UTS (MPa) | $\bar{\mu}_p$ ($\times 10^8$ Pa) | $\bar{\mu}_G$ ($\times 10^9$ Pa) | $\eta$ ($\times 10^{13}$ Pa.s) | Parameter loading value |
|---|---|---|---|---|---|
| 285°C | 43.0 | 4.79 | 1.43 | 3.95 | 0.45 UTS |
| 300°C | 40.2 | 4.12 | 0.51 | 2.23 | 0.45 UTS |
| 315°C | 36.3 | 4.19 | 0.79 | 4.04 | 0.30 UTS |
| 330°C | 23.8 | 5.07 | 0.79 | 3.19 | 0.20 UTS |

Table 1. Table showing values for the ultimate tensile strength (UTS) and various material parameters ($\bar{\mu}_p$, $\bar{\mu}_G$, $\eta$). The table also shows the loading value data set that was used to obtain the optimum set of material parameters.

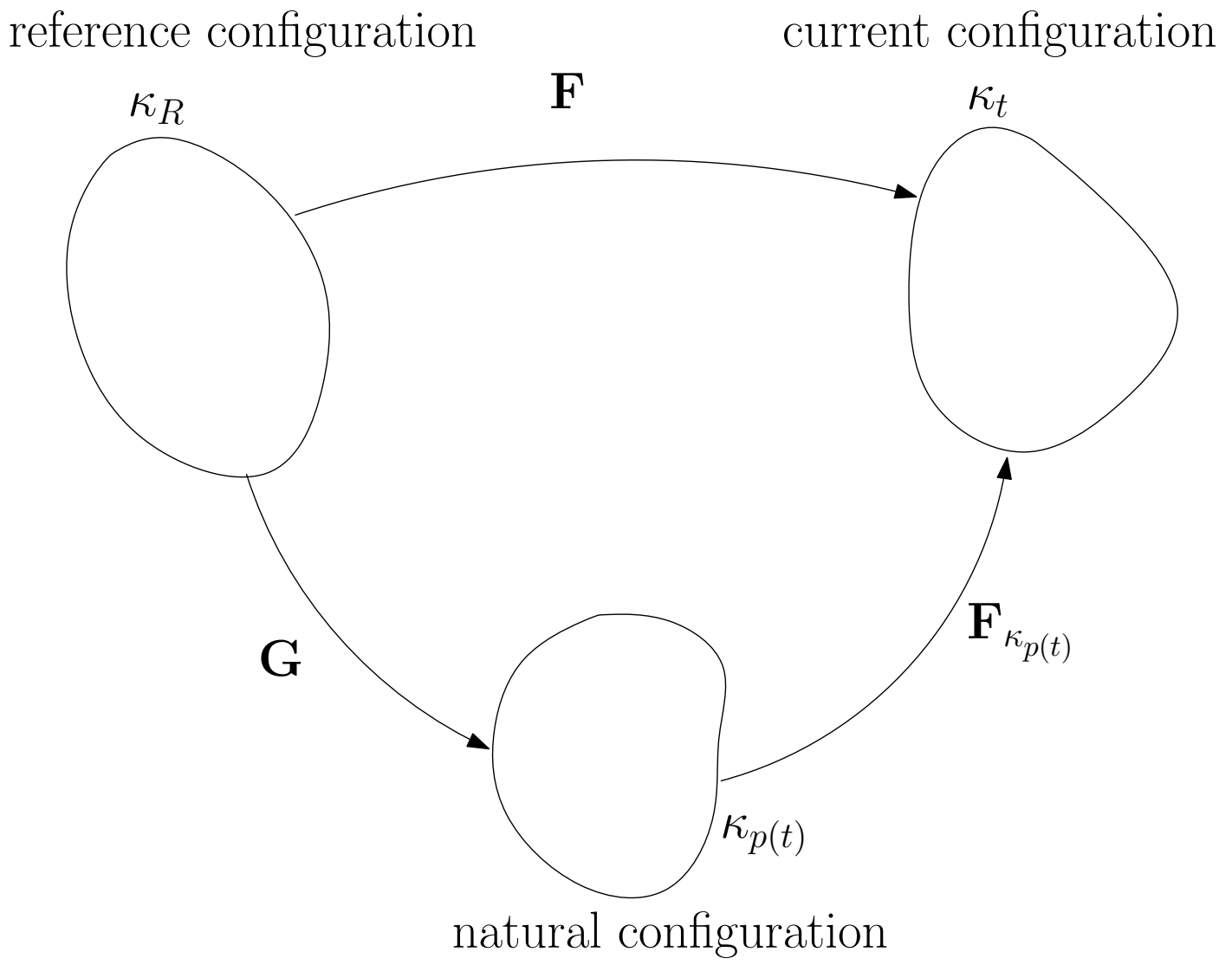

Figure 1. Illustration of various configurations of the body.

Satish Karra, Texas A&M University, Department of Mechanical Engineering, 3123 TAMU, College Station TX 77843-3123, United States of America
*E-mail address*: satkarra@tamu.edu

K. R. Rajagopal (Corresponding author), Texas A&M University, Department of Mechanical Engineering, 3123 TAMU, College Station TX 77843-3123, United States of America, Phone: 1-979-862-4552, Fax: 1-979-845-3081
*E-mail address*: krajagopal@tamu.edu

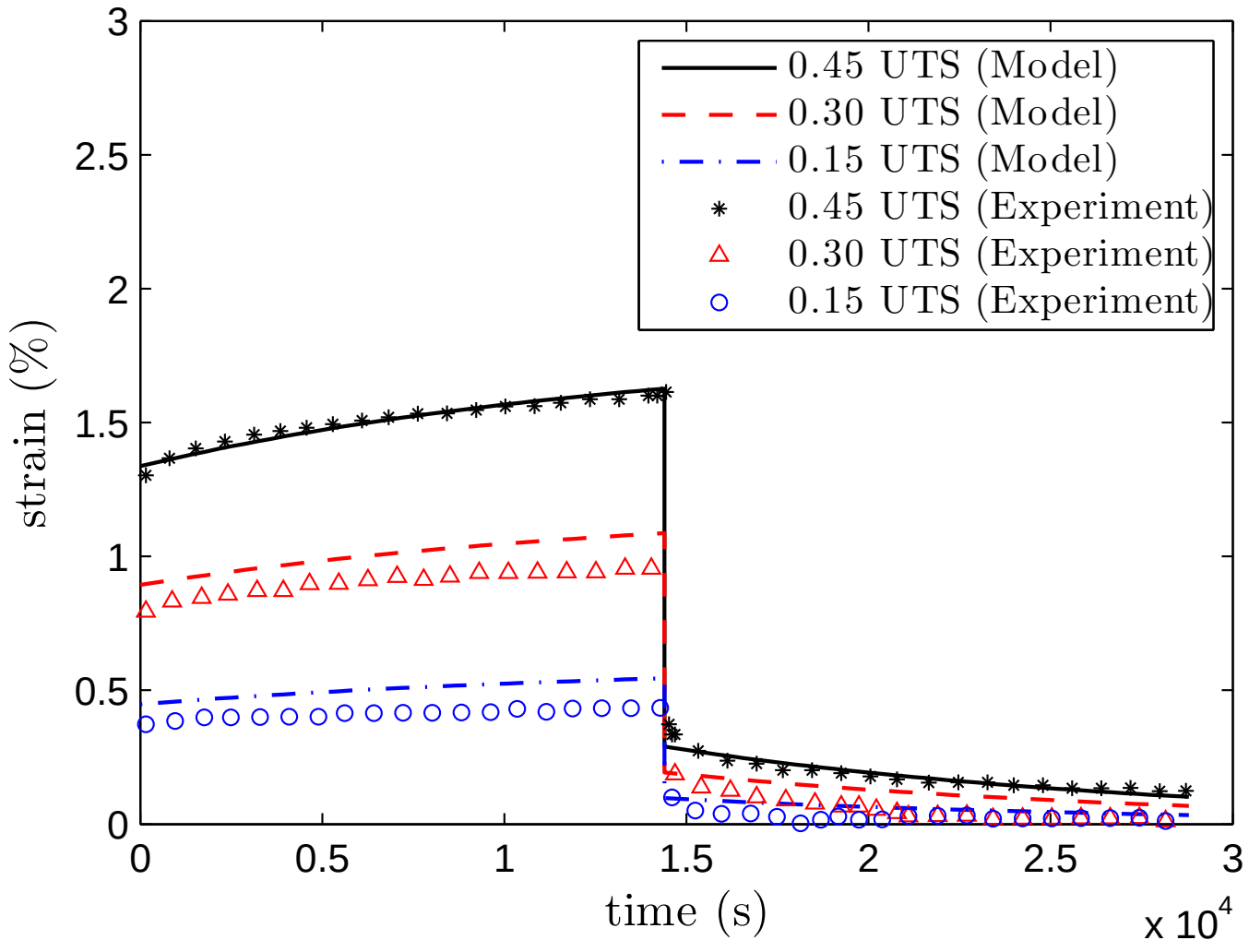

(A) 285°C

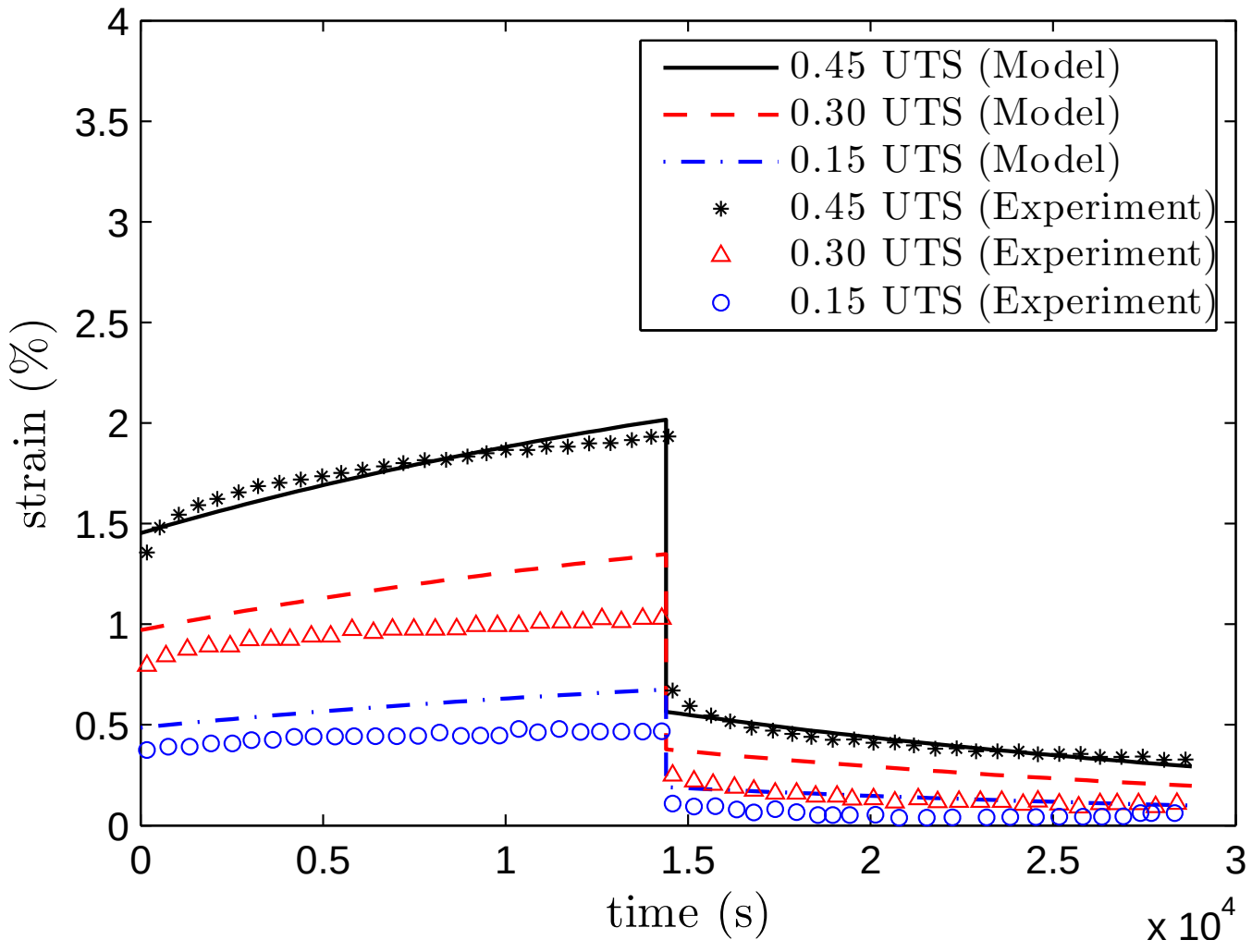

(B) 300°C

FIGURE 2. Comparison of the model predictions with experimental creep data of Bhargava (2007) at different loadings. The polyimide in this case is HFPE-II-52 at 285°C and 300°C. The parameters chosen and the values for the ultimate tensile strength (UTS) are shown in table (1). A weight of $w = 0.5$ was used for these two cases to obtain the optimum set of parameters.

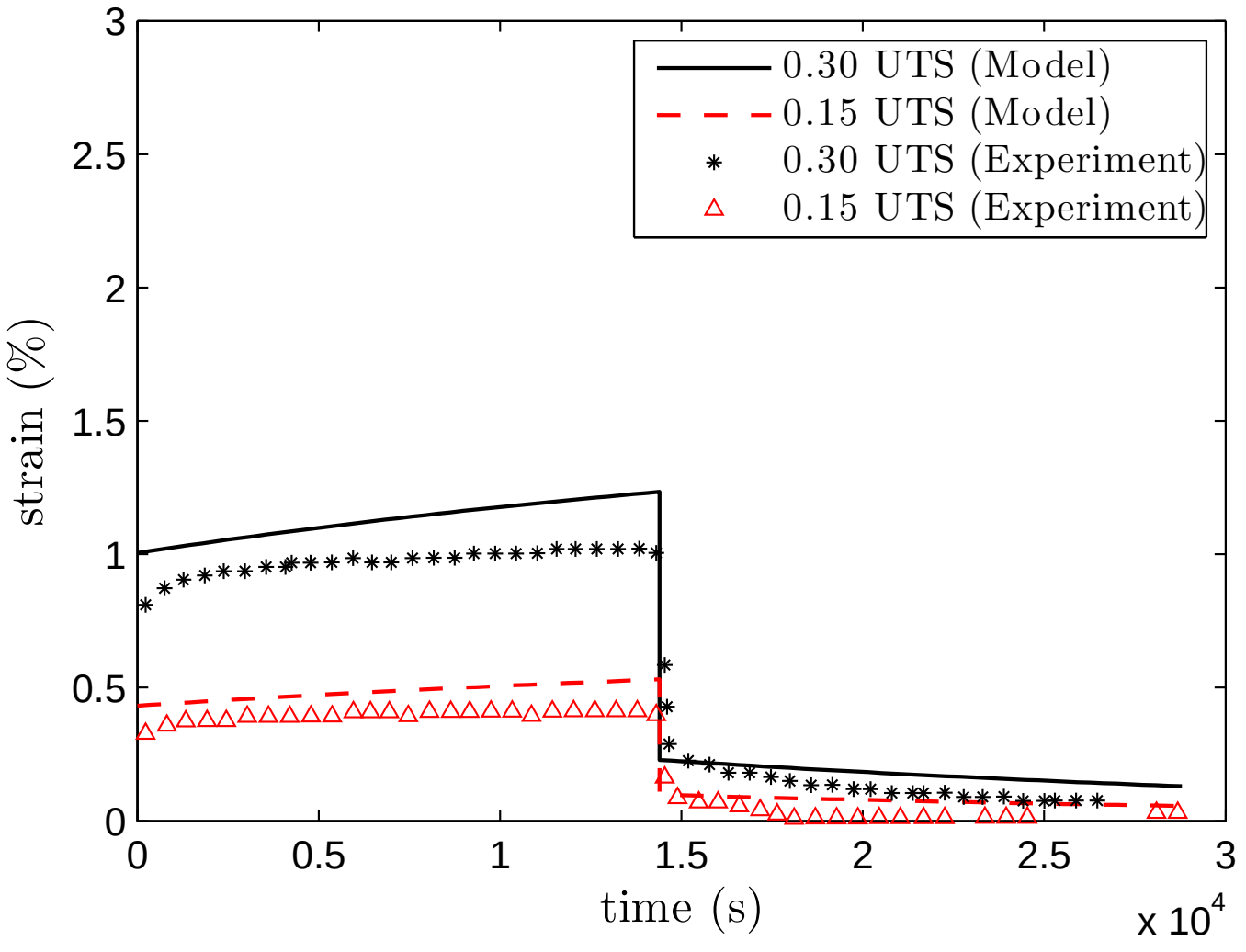

(A) 315°C

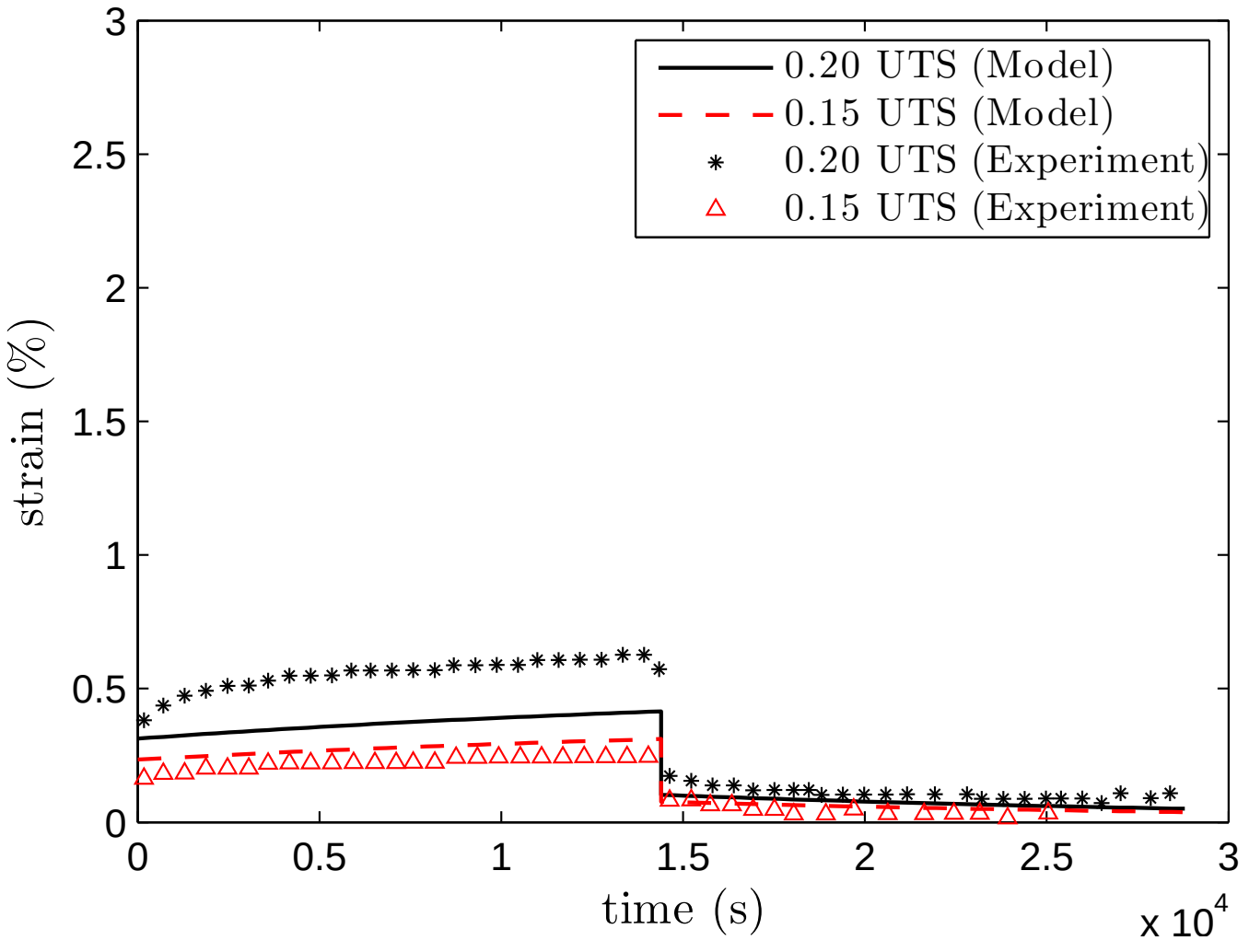

(B) 330°C

FIGURE 3. Comparison of the model predictions with experimental creep data of Bhargava (2007) at different loadings. The polyimide in this case is HFPE-II-52 at 315°C and 330°C. The parameters chosen and the values for the ultimate tensile strength (UTS) are shown in table (1). A weight of $w = 0.75$ was used for these two cases to obtain the optimum set of parameters.

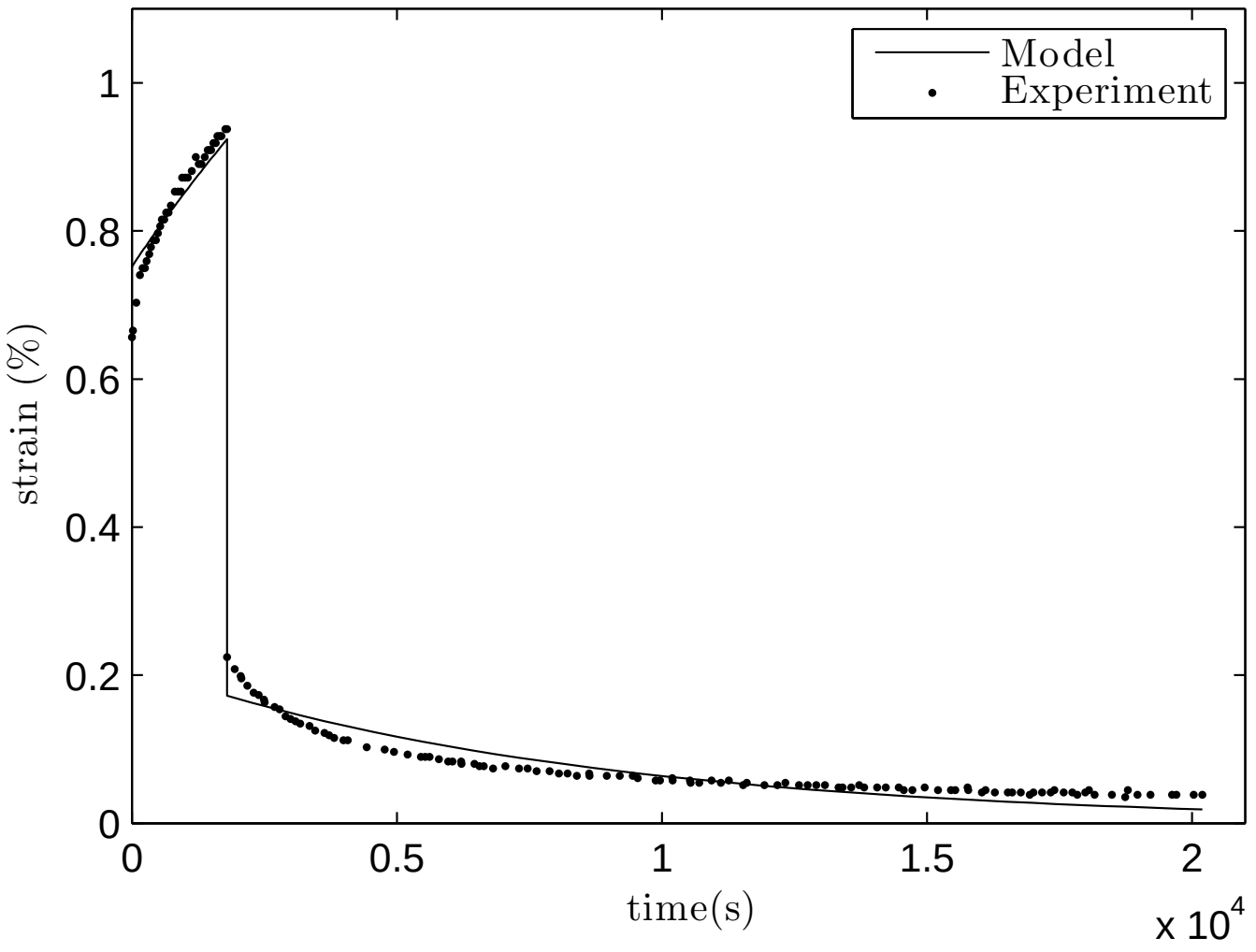

FIGURE 4. Comparison of the model with experimental creep data of Falcone and Ruggles-Wrenn (2009) for a loading of 10 MPa. The polyimide in this case is PMR-15 at a temperature of 288°C. The parameter values used were $\bar{\mu}_G = 4.42 \times 10^8$ Pa, $\bar{\mu}_p = 3.76 \times 10^8$ Pa, $\eta = 6.22 \times 10^{12}$ Pa.